\documentclass[12pt,twoside]{article}
\usepackage{amsmath,amssymb}
                                                     \newif\iffigs\figstrue
\iffigs
\input epsf.tex
\fi

\def\lskipamount{12pt}
\def\lskip{\vskip\lskipamount plus3pt minus2pt}
\def\lbreak{\par \ifdim\lastskip<\lskipamount
\removelastskip \penalty-200 \lskip \fi}

\def\lnobreak{\par \ifdim\lastskip<\lskipamount
\removelastskip \penalty200 \lskip \fi}

\makeatletter
\renewcommand{\@maketitle}{%
  \newpage
  \null
  \vskip 2em%
  \begin{center}%
    {\bf\Large \@title \par}%
    \vskip 1.5em%
    {\it\large
      \lineskip .5em%
      \begin{tabular}[t]{c}%
        \@author
      \end{tabular}\par}%
   \end{center}%
  \par
  \vskip 1.5em}
\makeatother

\makeatletter
\renewcommand{\@seccntformat}[1]{\csname the#1\endcsname. }
\renewcommand{\section}{\@startsection
    {section}{1}{0pt}%
    {-3.5ex \@plus -1ex \@minus -.2ex}%
    {2.3ex \@plus.2ex}%
    {\reset@font\normalsize\bfseries\centering}}
\renewcommand{\subsection}{\@startsection
    {subsection}{2}{0pt}%
    {-3.25ex \@plus -1ex \@minus -.2ex}%
    {1.5ex \@plus.2ex}%
    {\reset@font\normalsize\bfseries\centering}}
\renewcommand{\subsubsection}{\@startsection
    {subsubsection}{3}{0pt}%
    {-3.25ex \@plus -1ex \@minus -.2ex}%
    {1.5ex \@plus.2ex}%
    {\reset@font\normalsize\bfseries\centering}}
\renewcommand{\paragraph}{\@startsection
    {paragraph}{4}{0pt}%
    {3.25ex \@plus 1ex \@minus .2ex}%
    {-1em}%
    {\reset@font\normalsize\bfseries\centering}}
\renewcommand{\subparagraph}{\@startsection
    {subparagraph}{5}{0pt}%
    {3.25ex \@plus 1ex \@minus .2ex}%
    {-1em}%
    {\reset@font\normalsize\bfseries\centering}}
\makeatother

\def\thebibliography#1{\vskip 1.5pc{\centerline {\bf References}}\vskip 4pt
\list
 {[\arabic{enumi}]}{\settowidth\labelwidth{[#1]}\leftmargin\labelwidth
 \advance\leftmargin\labelsep
 \usecounter{enumi}}
 \def\newblock{\hskip .11em plus .33em minus .07em}
 \sloppy\clubpenalty4000\widowpenalty4000
 \sfcode`\.=1000\relax}

\def\affil#1{\vskip 1truepc{\noindent {\footnotesize #1}}}
\def\qed{\hbox to\textwidth{\hfill $\blacksquare$}\vskip 1.5truepc}

\makeatletter

\def\pplogo{\vbox{\kern-\headheight\kern-12pt
\halign{##&##\hfil\cr &\ppenumber\cr \rule{0pt}{2ex}&{\sc \ppnumber}\cr
\rule{0pt}{2.5ex}&\ppdate\cr}}}
\renewcommand{\maketitle}{\par
  \begingroup
    \renewcommand{\thefootnote}{\fnsymbol{footnote}}%
    \def\@makefnmark{\hbox to\z@{$\m@th^{\@thefnmark}$\hss}}%
    \long\def\@makefntext##1{\parindent 1em\noindent
            \hbox to1.8em{\hss$\m@th^{\@thefnmark}$}##1}%
    \if@twocolumn
      \ifnum \col@number=\@ne
        \@maketitle
      \else
        \twocolumn[\@maketitle]%
      \fi
    \else
      \newpage
      \global\@topnum\z@   
      \@maketitle
    \fi
    \thispagestyle{firstpage}\@thanks
  \endgroup
  \setcounter{footnote}{0}%
  \let\thanks\relax
  \let\maketitle\relax\let\@maketitle\relax
  \gdef\@thanks{}\gdef\@author{}\gdef\@title{}}

\def\ps@firstpage{\ps@plain \def\@oddhead{\hss\pplogo}%
  \let\@evenhead\@oddhead 
}

\makeatother

\def\ppnumber{duk-m-97-04}
\def\ppenumber{alg-geom/9705028}
\def\ppdate{May, 1997}

\begin{document}
\setcounter{page}{1}

\markboth%
{\it\normalsize David R. Morrison}%
{\it\normalsize Through the Looking Glass}

\title{Through the Looking Glass}
\author{David R. Morrison%
\thanks{Research partially supported by  National Science Foundation
grant DMS-9401447.}}
\maketitle

\newtheorem{conjecture}{Conjecture}\renewcommand\theconjecture{}

\newcommand\SU{\operatorname{SU}}
\newcommand\Aut{\operatorname{Aut}}
\newcommand\Hom{\operatorname{Hom}}
\newcommand\Pic{\operatorname{Pic}}
\newcommand\Hull{\operatorname{Hull}}
\newcommand\Span{\operatorname{Span}}
\renewcommand\Im{\operatorname{Im}}

\renewcommand\wp{{\Bbb P}^{(1,1,2,2,2)}}
\newcommand\p{{\Bbb P}}

\newcommand{\gitquot}{\mathchoice
{\mathrel{\mskip-4.5mu/\!/\mskip-4.5mu}}
{\mathrel{\mskip-4.5mu/\!/\mskip-4.5mu}}
{\mathrel{\mskip-3mu/\mskip-4.5mu/\mskip-3mu}}
{\mathrel{\mskip-3mu/\mskip-4.5mu/\mskip-3mu}}
}

\begin{abstract}

It is frequently possible to produce new Calabi--Yau threefolds from
old ones by a process of allowing the complex structure to degenerate
to a singular one,
and then performing a resolution of singularities.
(Some care is needed to ensure that the Calabi--Yau condition be preserved.)
There has been speculation that all Calabi--Yau threefolds could be
linked in this way, and considerable evidence has been amassed in this
direction.
We propose here a natural way to relate this construction
to the string-theoretic phenomenon known as ``mirror symmetry.''  We
formulate a
conjecture which in principle could predict mirror partners
for all Calabi--Yau threefolds, provided that all were indeed linked
by the degeneration/resolution process.  The conjecture produces
new mirrors from old, and so requires some initial mirror manifold
construction---such as Greene-Plesser orbifolding---as a starting point.
(Lecture given at the CIRM conference, Trento, June 1994, and at the Workshop
on Complex Geometry and Mirror Symmetry, Montr\'eal, March 1995.)
\end{abstract}


\section{Introduction} \label{sec:1}

One of the most intriguing conjectures in complex analytic geometry
is the spe\-cu\-la\-tion---essentially due to Herb Clemens, with
refinements by Miles Reid and others---that
the moduli spaces of all nonsingular compact complex threefolds
with trivial canonical bundle might possibly
be connected into a single family.  Clemens introduced\footnote{The
construction is hinted at in  \cite{clmnone, clmntwo}, and given
in detail in \cite{friedman}; the consequences for connecting moduli
spaces are spelled out in \cite{reid} and \cite{friedmansurvey}, and the
local geometry of the transition was analyzed in \cite{cd}.}
a transition process among the moduli spaces of these various threefolds
in which the complex structure on one such threefold
is allowed to acquire singularities,
and  those singularities are then resolved by means of
a bimeromorphic map from a second such threefold.  If the threefolds involved
have compatible K\"ahler metrics (which therefore determine
``Calabi--Yau structures''), then the second
step in this process can be viewed as a sort of inverse to
allowing the K\"ahler metric on the
nonsingular model to degenerate to some kind of
``K\"ahler metric with singularities,'' the singularities being concentrated
on the exceptional set of the resolution map.

It is natural to wonder how this process might be related to the
string-theoretic phenomenon
known as ``mirror symmetry'' \cite{dixon, lvw, cls, gp},
which predicts that the moduli space of complex structures on one
Calabi--Yau threefold $X$ should in many cases be locally isomorphic to
the space of (complexified) K\"ahler structures on a ``mirror
partner'' $Y$ of $X$, and {\it vice versa}.
Viewed through the ``looking glass'' of mirror symmetry,
the two ways of approaching the singular space in Clemens' transition appear
complementary---one involves a specialization of complex structure parameters,
the other involves a specialization of K\"ahler parameters---and
indeed their r\^oles should be reversed when passing to mirror partners.

The thesis of this lecture is that the whenever such a transition exists, it
ought to enable us to predict the mirror partner of one of the Calabi--Yau
manifolds involved in the transition from a knowledge of the
mirror partner of the other.  We will explain this idea in general terms,
give some evidence and examples, and then formulate a specific
conjecture which implements it.

After this lecture was delivered, there were a number of
new developments in string theory related to this construction; we describe
those briefly at the end.

\section{Extremal Transitions and Mirror Symmetry }\label{sec:2}

Use of the term
``Calabi--Yau'' varies; our conventions are as follows.  Let
$X$ be a compact, connected, oriented manifold of dimension $2n$.
A {\it Calabi--Yau metric}\/ on $X$ is a Riemannian metric
whose (global)
holonomy is a subgroup of $\SU(n)$.  A {\it Calabi--Yau structure}\/
on $X$
is a choice of complex structure with trivial canonical bundle
together with
a K\"ahler metric; each Calabi--Yau structure determines
a unique Calabi--Yau metric
according to Yau's solution \cite{Yau} of the Calabi conjecture
\cite{Calabi}.  Finally, $X$ is
a {\it Calabi--Yau manifold}\/ if it admits
Calabi--Yau structures (and hence Calabi--Yau metrics).

Let $X$ be a Calabi--Yau manifold. The simplest way to specify a family
of Calabi--Yau structures on $X$ is by means of a proper smooth
holomorphic map $\pi:{\cal X}\to S$ satisfying $\omega_{{\cal X}/S}={\cal
O}_S$ whose fibers
${\cal X}_s:=\pi^{-1}(s)$ are diffeomorphic to $X$, together with
a fixed K\"ahler metric on the total space ${\cal X}$,
or more generally, a family of K\"ahler metrics depending on parameters.
We can allow
the complex structures to acquire singularities by enlarging the
family to $\overline\pi:\overline{\cal X}\to\overline{S}$,
where $S\subset\overline{S}$ is an open subset whose complement is closed, and
where
$\overline\pi$ is still assumed to be a proper holomorphic map such that
$\omega_{\smash{\overline{\cal X}/\overline S}}=
{\cal O}_{\smash{\overline S}}$,
but
is no longer assumed to be smooth.
The fibers   $\overline{\cal X}_{\sigma}:=\overline\pi^{-1}(\sigma)$ for
$\sigma\in \widehat{S}:=\overline{S}-S$
may then be complex analytic spaces
with singularities; of course, such spaces do not have K\"ahler metrics
{\it per se}.

Suppose that the singularities of the fibers $\overline{\cal X}_{\sigma}$
can be simultaneously resolved for $\sigma\in \widehat{S}$ by
manifolds with Calabi--Yau structures.\footnote{A familiar instance of
this is the small resolution of a collection of ordinary double points
\cite{Atiyah, Hirz, TianYau}; note that we are demanding that a K\"ahler
metric exist
on the resolved space, which puts restrictions on the
global configuration of the singularities.}  That is,
suppose that there exists
a manifold $\widehat X$, a family of complex structures
$\widehat\pi:\widehat{\cal X}\to \widehat{S}$ on $\widehat{X}$ such
that $\omega_{\smash{\widehat{\cal X}/\widehat{S}}}={\cal
O}_{\smash{\widehat{S}}}$ whose total space $\widehat{\cal X}$ is a
K\"ahler manifold, and a proper bimeromorphic map $f:\widehat{\cal X}\to
\overline{\cal X}_{\smash{\widehat{S}}}
:=\overline\pi^{-1}(\widehat{S})$ which commutes with projection to
$\widehat{S}$.  (In general, it will be necessary to shrink
$\overline{S}$ before this simultaneous resolution is possible---if it is
possible at all.)  If in addition, the bimeromorphic map
$f_{\sigma}:\widehat{\cal X}_{\sigma}\to\overline{\cal X}_{\sigma}$ is an
{\it extremal contraction}\/ in the sense of Mori theory (which means that
for every
algebraic curve $C\subset \widehat{\cal X}_{\sigma}$, we have $\dim f(C)=0$
if and only if the class of $C$ lies in some fixed face ${\cal F}$ of the
Mori cone of $\widehat{\cal X}_{\sigma}$), then we say that $X$ and
$\widehat{X}$ are related by an {\it extremal transition}.

We can also describe an extremal transition ``in reverse,'' starting with
$\widehat{X}$.  In this version of the definition, we should begin with
a family $\widehat{\cal X}\to \widehat{S}$ of complex structures on
$\widehat{X}$
whose total space is a K\"ahler manifold, fix a birational extremal
contraction $\widehat{\cal X}_\sigma\to \overline{\cal X}_\sigma$ (which is
determined simply by the choice of an appropriate face ${\cal F}$ of the
Mori cone), and ask for a smoothing of the contracted spaces.  That is, we
look for a space $\overline{S}$
into which $\widehat{S}$ can be embedded as a closed subset, and a family
$\overline{\cal X}\to \overline{S}$ such that $\overline{\cal X}_s$ is
smooth for $s\in \overline{S}-\widehat{S}$ (and diffeomorphic to $X$)
whereas $\overline{\cal
X}_{\smash{\widehat{S}}}$ is the same family of singular spaces as
before. If this smoothing
exists, then $X$ and $\widehat{X}$ are related by an extremal transition as
above.

This ``reverse'' description can actually be thought of as varying the
K\"ahler parameters, rather than the complex structure parameters.
Consider a family of K\"ahler metrics on the total space $\widehat{X}$ whose
cohomology class approaches a point on the dual ${\cal F}^\perp$ of the
face ${\cal F}$.  Metrically, the curves $C$ lying in that face will
approach zero area in such a process, and so the space with metric included
``approaches'' the contracted space $\overline{\cal X}$.

In string theory, the choice of K\"ahler class serves as an important
parameter in the theory.  In fact, a more natural parameter for string theory
is a complexification of the
K\"ahler class, varying throughout an open subset of
the {\it complexified K\"ahler moduli space}\/
${\cal M}_{\text{K\"ah}}(\widehat{\cal X}_\sigma):= {\cal K}_{\Bbb
C}(\widehat{\cal X}_\sigma)/\Aut(\widehat{\cal X}_\sigma)$,
where
\begin{equation}
{\cal K}_{\Bbb C}(\widehat{\cal X}_\sigma):=
\{B+i\omega\in H^2(\widehat{X},{\Bbb C}/{\Bbb Z})\ |\ \text{$\omega$ is a
K\"ahler class on $\widehat{\cal X}_\sigma$}\}
\end{equation}
is the {\it complexified K\"ahler cone}, and $\Aut(\widehat{\cal
X}_\sigma)$ is the group of holomorphic automorphisms.
The face ${\cal F}$ of the Mori cone corresponding to an
extremal contraction determines a boundary ``wall''
${\cal F}^{\perp}$ of the complexified K\"ahler cone,
and---provided that the singularities of $\overline{{\cal X}}_\sigma$ are
sufficiently mild---the intersection
\begin{equation}
{\cal F}^\perp\cap\overline{{\cal K}_{\Bbb C}(\widehat{\cal X}_\sigma)}
\end{equation}
will coincide with
$\overline{{\cal K}_{\Bbb C}({\cal X}_s)}$
for generic nearby $s\in S$.
If the automorphism groups are reasonably well-behaved, there is then a
natural inclusion
\begin{equation}\label{inclusion}
\overline{{\cal M}_{\text{K\"ah}}({\cal X}_s)}\subset
\overline{{\cal M}_{\text{K\"ah}}(\widehat{\cal X}_\sigma)}
\end{equation}
(using the partial compactifications described in \cite{compact}).

{\it Mirror symmetry}\/ refers to a phenomenon in string theory
in which certain pairs of
Calabi--Yau manifolds produce isomorphic physical theories, in such a way
that the complexified K\"ahler moduli space of the first Calabi--Yau
manifold is mapped to the ordinary (complex structure) moduli space of the
second and {\it vice versa}, establishing local isomorphisms between
those moduli spaces.
Our basic proposal is that
the mirror of an extremal transition should be another such transition, from
a smooth space $Y$ specializing to a singular space $\overline{Y}$ (using
the complex structure) which has a resolution of singularities $\widehat{Y}$.
In this mirror version,
{\it the mirror partner of $X$ should be $\widehat{Y}$, while the mirror
partner of $\widehat{X}$ should be $Y$}.
Under the mirror map, the inclusion
\begin{equation}\label{inclX}
\widehat{S}\subset\overline S
\end{equation}
between compactified parameter spaces for complex structures on
$\widehat{X}$ and $X$ should locally map to the inclusion
\begin{equation}\label{inclY}
\overline{{\cal M}_{\text{K\"ah}}({\cal Y}_t)}\subset
\overline{{\cal M}_{\text{K\"ah}}(\widehat{\cal Y}_\tau)}
\end{equation}
(and similarly for complex structures on $\widehat{Y}$ and $Y$ mapping to
K\"ahler structures on $X$ and $\widehat{X}$).
In particular, if we know the mirror partner of
$X$ (resp.~$\widehat{X}$), and if we can predict where in the moduli space
the mirror of the extremal transition should occur, then we should be able to
determine the mirror partner of $\widehat{X}$ (resp.~$X$).
We will make this proposal more precise in section \ref{sec:4}.

\section{Evidence and Examples}\label{sec:3}

\subsection{ Hodge numbers of conifold transitions}\label{subsec:31}

In addition to the formal analogy between K\"ahler and complex structures
which led to our ``looking glass'' interpretation of the mirror of
an extremal transition,
we were motivated by a computation of the effect on cohomology
of the simplest type of extremal transition which has
been known since the work of Clemens. (See \cite{werner, schoen, WvG}
for general versions of this computation).
Suppose that the complex dimension of our Calabi--Yau manifolds is three,
and suppose that all singularities of
$\overline{\cal X}_v$ are ordinary double points.
(In this case, we will follow the conventions of the physics literature
\cite{ghone} and refer to the extremal transition as a
{\it conifold transition.}\footnote{We are assuming that the physical model
does not
have ``discrete torsion'' at the ordinary double points; otherwise, the
picture is somewhat different \cite{VW, stabs}.})
Generally, we might expect
that acquiring a double point places one condition on moduli.  However, the
double points may fail to impose independent conditions.  Clemens'
computation relates the failure of double points to
impose independent conditions on the moduli of $X$
 to the relative Picard number
of the small resolution $\widehat{X}$, in the following way:
if there are $\delta$ double points which only impose $\sigma:=\delta-\rho$
conditions on moduli then

\medskip

\noindent (i)
 $H^3(\widehat{X})\subset H^3(X)$ is a subspace of codimension $2\sigma$
(the inclusion comes from the identification of $H^3(\widehat{X})$ with the
space of invariant cycles), and

\medskip

\noindent (ii)
$\oplus H^{2k}(X) \subset \oplus H^{2k}(\widehat{X})$ is a subspace of
codimension $2\rho$ (the inclusion is induced by pullbacks of divisors).

\medskip

\noindent
In other words, each failure to impose a condition on moduli leads to
an extra class in $H^2(\widehat{X})$.

One of the properties of mirror symmetry for threefolds is that $H^3$ and
$\oplus H^{2k}$ are exchanged when passing to a mirror partner.  Thus,
if mirror partners $\widehat{Y}$ and $Y$ are known for $X$ and
$\widehat{X}$, respectively, then by using these mirror isomorphisms of
cohomology we find that

\medskip

\noindent (i)
$\oplus H^{2k}(Y) \subset \oplus H^{2k}(\widehat{Y})$ has codimension
$2\sigma$,
and

\medskip

\noindent (ii)
$H^3(\widehat{Y})\subset H^3(Y)$ has codimension $2\rho$.

\medskip

\noindent
 It then seems very natural to conjecture
that $Y$ and $\widehat{Y}$ should be related by a conifold transition
in the opposite direction:  it should be possible to allow the
complex structure on $Y$ to acquire $\delta$ ordinary double
points which this time
impose only $\rho$ conditions on moduli, such that the singular spaces
obtained can be resolved by $\widehat{Y}$.\footnote{Recent results
\cite{KMP} indicate that this conjecture can only hold for generic
moduli, if we insist on using only conifold transitions rather than
the more general extremal transitions.}
This would be a special case of our general principle
relating mirror symmetry to extremal transitions.

\subsection{Extremal transitions among toric hypersurfaces}\label{subsec:32}

One context in which certain extremal transitions have a very explicit
realization is the case of Calabi--Yau hypersurfaces in toric varieties.
An {\it $n{+}1$-dimensional toric variety}\/ $V$ is a $T$-equivariant
compactification of the algebraic torus $T:=({\Bbb C}^*)^{n+1}$ (often
specified by the combinatorial data encoded in a {\it fan}).  To describe a
hypersurface within such a variety we need a polynomial---the defining
equation---whose constituent monomials can be thought of as ${\Bbb
C}^*$-valued characters $\chi:T\to{\Bbb C}^*$.  Batyrev \cite{batyrev} has
given an elegant condition which characterizes when such a hypersurface
will be  Calabi--Yau.  The condition is stated in terms of the
{\it Newton polyhedron}\/ of the polynomial, which is the polyhedron
spanned by the monomials appearing in the equation.  That is, if the
polynomial is $\sum_{a\in M} c_a \chi_a$ where $M$ is the lattice of
characters on $T$, then the Newton polyhedron is the convex hull
\begin{equation}
{\cal P}:=\Hull(\{a\ |\ c_a\ne0\})\subset M\otimes{\Bbb R}.
\end{equation}
Batyrev's
criterion says that the generic such hypersurface is
Calabi--Yau\footnote{We must use a slight generalization of the term
``Calabi--Yau'' here, in which Gorenstein canonical singularities are
allowed on these hypersurface.  However, when $n\le3$ it is always possible to
choose the fan in such a way that the generic such hypersurface is
nonsingular, and we can then speak of Calabi--Yau {\it manifolds}\/ in this
context.} if ${\cal P}$ is
{\it reflexive}, which means that (1) ${\cal P}$ is the convex hull of the
lattice points it contains, (2) there is a unique lattice point $a_0$ in
the interior of ${\cal P}$, and (3) the polar polyhedron (with respect to
$a_0$) defined by
\begin{equation}
{\cal P}^\circ:=\{x\in N\otimes{\Bbb R}\ |\ \langle x,a\rangle - \langle
x,a_0 \rangle \ge -1\ \hbox{for all}\ a\in{\cal P}\}
\end{equation}
has its vertices at lattice points of the dual lattice $N:=\Hom(M,{\Bbb
Z})$.

Conversely, given a reflexive polyhedron ${\cal P}$, there is an associated
family of Calabi--Yau hypersurfaces with defining equation
\begin{equation}\label{poly}
f_{\cal P}=\sum_{a\in{\cal P}\cap M}c_a\chi_a
\end{equation}
embedded in a toric variety $V_{\cal P}$ determined by some fan whose
one-dimensional cones
are the rays ${\Bbb R}\vec{v}$ for
$\vec{v}\in {\cal P}^\circ\cap N$.
(The toric variety
$V_{\cal P}$ is not uniquely specified by this,
since we need to choose the fan, not just the one-dimensional
cones; however, $\Pic(V_{\cal P})$ is independent of the choice of fan.)

There is a simple class of {\it toric extremal transitions}\/ which can be
described in
these terms.\footnote{This construction was independently noted by
Berglund, Katz and Klemm \cite{bkk}.}  Suppose we have two reflexive
polyhedra ${\cal Q}\subset{\cal P}$.  Since all monomials appearing in
$f_{\cal Q}$ also appear in $f_{\cal P}$, we can regard the hypersurfaces
$X_{\cal Q}\subset V_{\cal Q}$
associated to ${\cal Q}$ as being limits of the hypersurfaces $X_{\cal
P}\subset V_{\cal P}$ associated to
${\cal P}$.  In fact, they will have worse than generic singularities when
embedded in
the toric variety $V_{\cal P}$, but those singularities can be improved, or
in some cases (including the case $n\le3$) fully resolved, by
further triangulation.  In fact, the unique interior point $a_0$ of
${\cal Q}$ must also be the unique interior point of ${\cal P}$, which
implies that
\begin{equation}
{\cal P}^\circ \subset {\cal Q}^\circ.
\end{equation}
Starting with a triangulation of
${\cal P}^\circ$, we will be able to further triangulate by including
vertices of ${\cal Q}^\circ$.

This construction provides an extremal transition as defined above, when
the spaces involved are nonsingular (as is the case when $n\le3$): the
family of complex structures acquires canonical singularities when certain
coefficients in the defining equation are set to zero, and those
singularities can be simultaneously resolved by maps which are extremal
contractions.  More precisely, the complex structures on the Calabi--Yau
manifold $X_{\cal P}$ have a
natural parameter space $S$ which is an open subset of
\begin{equation}\label{barS}
\overline{S}:={\Bbb C}^{{\cal P}\cap M} \gitquot \Aut(V_{\cal
P})\times{\Bbb C}^*.
\end{equation}
(The notation indicates a quotient in the sense of Geometric Invariant
Theory.)
One of the subsets of $\overline{S}$ along which the Calabi--Yau spaces
approach a singular space $\overline{X}$ is the set
\begin{equation}\label{hatS}
\widehat{S}:={\Bbb C}^{{\cal Q}\cap M} \gitquot \Aut(V_{\cal Q})\times{\Bbb
C}^*
\end{equation}
in which all coefficients $c_a$ in eq.~\eqref{poly} with $a\not\in{\cal Q}$
have been set to zero.  This is where our extremal transition is located.
Refining a triangulation of ${\cal P}^\circ$ to a triangulation of ${\cal
Q}^\circ$ produces the resolution of singularities
$\widehat{X_{\cal Q}}\to\overline{X}$.

In addition to giving a criterion for when a toric hypersurface is
Calabi--Yau, Batyrev made a simple and beautiful
conjecture concerning a possible mirror partner for any
such hypersurface.  To find the mirror family for $X_{\cal P}$, Batyrev's
conjecture states that we should simply use
the polar polyhedron ${\cal P}^\circ$ to determine a new family of
hypersurfaces---let's call it $\widehat{Y}_{{\cal P}^\circ}$.  Batyrev's
construction exhibits a perfect compatibility with toric extremal
transitions:
since ${\cal P}^\circ \subset {\cal Q}^\circ$, we {\it automatically}\/
get an extremal transition between $Y_{{\cal Q}^\circ}$
and $\widehat{Y}_{{\cal P}^\circ}$.  Moreover, Batyrev's formula
\cite{batyrev} for the Hodge numbers of these hypersurfaces
shows that---as in the case of our conifold
transition conjecture---mirror symmetry is compatible with extremal
transitions in the manner stated.
In fact, using a description of the K\"ahler cones of $Y_{{\cal Q}^\circ}$ and
$\widehat{Y}_{{\cal P}^\circ}$ in terms of ${\cal Q}$ and ${\cal P}$, together
with an extension of Batyrev's conjecture known as the ``monomial-divisor
mirror map'' \cite{mondiv}, one can verify that the inclusion
$\widehat{S}\subset\overline{S}$ between the spaces from eqs.~\eqref{hatS}
and \eqref{barS} is locally isomorphic to the inclusion\footnote{Actually,
we are only working with a subset of the K\"ahler moduli space here
corresponding to ``toric'' divisors---divisors which arise by restriction
from a divisor on the ambient toric variety.  Similarly, the parameter
spaces $\overline{S}$ and $\widehat{S}$ only capture those complex
structure parameters which preserve the property that the Calabi--Yau space
can be embedded as a hypersurface in a toric variety.}
\begin{equation}
\overline{{\cal M}_{\text{K\"ah}}(Y_{{\cal Q}^\circ})}\subset
\overline{{\cal M}_{\text{K\"ah}}(\widehat{Y}_{{\cal P}^\circ})}
\end{equation}
between complexified K\"ahler moduli spaces, with the corresponding face
${\cal F}^\perp$ of the K\"ahler cone determined by the inclusion ${\cal
Q}\subset{\cal P}$.

\subsection{An explicit example}\label{subsec:33}

The relationship between mirror symmetry and extremal transitions can
be seen very clearly in an explicit example worked out some years ago
by the author in collaboration with
P.~Candelas, X. de~la~Ossa, A.~Font and S.~Katz
\cite{cdfkm}.\footnote{We shall use notation compatible with
\cite{Small, MP}, where further details about this example were analyzed.}
The transition is most easily described from the ``reverse'' perspective,
beginning with $\widehat{Y}$ and varying the K\"ahler class to eventually
produce $Y$.

We begin with the weighted projective space $\wp$, and let
$\pi: \widehat{\p}\to\wp$
be the blowup of the singular locus of $\wp$.  Then $\widehat{\p}$ is a smooth
fourfold containing smooth anti-canonical divisors $\widehat{Y}$, which are
Calabi--Yau threefolds.  The induced blowdown $\pi:\widehat{Y}\to\overline{Y}$
is an extremal contraction.  Moreover, the Calabi--Yau space
$\overline{Y}$ (which is defined by a homogeneous polynomial of weighted
degree $8$, and has canonical singularities) can be smoothed.  This is most
easily seen by using
the linear system ${\cal O}(2)$ to map $\wp$ to $\p^5$, where the
image is a quadric hypersurface of rank 3.  The image of $\overline{Y}$
under this mapping is the intersection of this singular quadric hypersurface
with a smooth hypersurface of degree $4$.  This has a smoothing to
a space $Y$, the intersection of smooth hypersurfaces of degrees $2$ and $4$.

To describe $\widehat{\p}$ in standard toric geometry language, we begin
with the description of the weighted projective space $\wp$, which is
determined by the following lattice points in $N={\Bbb Z}^4$:
\begin{equation*}
\begin{array}{ll}
\multicolumn{2}{c}{v_1=(-1,-2,-2,-2),}\\
v_2=(1,0,0,0),&v_4=(0,0,1,0),\\
v_3=(0,1,0,0),&v_5=(0,0,0,1).
\end{array}
\end{equation*}
The blowup is described by including the additional lattice point
\begin{equation*}
v_6=\frac{v_1+v_2}2=(0,-1,-1,-1).
\end{equation*}
Each lattice point $v_i$ determines a toric divisor $D_i$, whose classes
generate $\Pic(\widehat{\p})$.  If we choose a basis for the linear
relations among the $v_i$'s, say
\begin{align*}
v_3+v_4+v_5+v_6 & =0\\
v_1+v_2-2v_6 & =0
\end{align*}
then there is a corresponding basis $\eta_1$, $\eta_2$ of $\Pic(\widehat{\p})$
for which
\begin{align*}
[D_3]=[D_4]=[D_5]&=\eta_1\\
[D_1]=[D_2]&=\eta_2\\
[D_6]&=\eta_1-2\eta_2.
\end{align*}
It turns out that $\eta_1$ and $\eta_2$ generate the K\"ahler cone of
$\widehat{\p}$.

We can write an arbitrary K\"ahler class in the form
$a_1\eta_1+a_2\eta_2$, and use $a_1$, $a_2$ as coordinates on the K\"ahler
cone; the natural coordinates to use on the complexified K\"ahler cone are
then $\exp(2\pi i\,a_1)$, $\exp(2\pi i\, a_2)$.  The K\"ahler cone has two
faces given by $\{a_1=0\}$ and $\{a_2=0\}$.  The first face corresponds to a
pencil of K3 surfaces, and has no associated extremal transition.  It is
the second face $\{a_2=0\}$ which is associated to our extremal transition.
The class $\eta_1$ which spans that face
contains $2D_1+D_6, D_1+D_2+D_6,2D_2+D_6,D_3, D_4, D_5$
as representatives, and the corresponding linear system maps $\widehat{\p}$
to $\wp$, shrinking
the exceptional
divisor of the blowup map to zero size. As remarked above, the
hypersurfaces $\overline{Y}$
in $\wp$ can then be smoothed (after reembedding $\wp$ as a rank 3 quadric in
$\p^5$); this gives a complete description of the extremal transition.

Candidate mirror partners are known for both $\widehat{Y}$ and $Y$
\cite{gp, lt}, so
it is natural to look for a connection between these---a mirror image of
the extremal transition mentioned above.  This was also found in \cite{cdfkm}.
The candidate mirror partner of $\widehat{Y}$ is the desingularization of
an anti-canonical hypersurface in $\wp/G$, where $G$ is the image in
$\Aut(\wp)$ of
\begin{equation*}
\widetilde{G}=\{
\vec{\lambda}
\in({\Bbb C}^*)^5\ |\
\lambda_1^8=\lambda_2^8=\lambda_3^4=\lambda_4^4=\lambda_5^4=
\lambda_1\lambda_2\lambda_3\lambda_4\lambda_5=1\}.
\end{equation*}

The complex structure moduli space of the mirror manifold $X$ can be described
in terms of some analogous parameters.  Rather than using a redundant
description
in terms of toric divisors, we can this time use a redundant description
in terms of monomials.  To describe the mirror of $\widehat{Y}$, in fact,
we should consider the family of polynomials in the homogeneous coordinates
$x_1$, \dots, $x_5$
\begin{equation}
c_0 x_1x_2x_3x_4x_5 + c_1x_1^8 + c_2x_2^8 + c_3x_3^4 + c_4x_4^4 + c_5x_5^4
+ c_6 x_1^4x_2^4.
\end{equation}
The connection with the divisors on the mirror becomes more apparent if we
divide the polynomial by $x_1x_2x_3x_4x_5$ and rewrite in terms of the
basis of the torus $T=({\Bbb C}^*)^4$ defined by
\begin{align*}
t_1&=x_1^{-1}x_2^{7}x_3^{-1}x_4^{-1}x_5^{-1}\\
t_2&=x_1^{-1}x_2^{-1}x_3^{3}x_4^{-1}x_5^{-1}\\
t_3&=x_1^{-1}x_2^{-1}x_3^{-1}x_4^{3}x_5^{-1}\\
t_4&=x_1^{-1}x_2^{-1}x_3^{-1}x_4^{-1}x_5^{3}.
\end{align*}
In this basis, the defining polynomial of $X$ becomes
\begin{equation}
c_0+c_1t_1^{-1}t_2^{-2}t_3^{-2}t_4^{-2}
+c_2t_1+c_3t_2+c_4t_3+c_5t_4+c_6t_2^{-1}t_3^{-1}t_4^{-1}
\end{equation}
 and the exponents on the monomials correspond to the $v_i$'s of the mirror.
The coordinates on the complex structure moduli space which are analogous
to the $\exp(2\pi i\,a_j)$'s are then
\begin{equation}
q_1:=c_3c_4c_5c_6/c_0^4;\qquad q_2:=c_1c_2/c_6^2.
\end{equation}
(These are local coordinates on the parameter space
$\overline{S}={\Bbb C}^{7} \gitquot T\times{\Bbb C}^*$.)

\iffigs
\begin{figure}[t]
\begin{center}
\begin{picture}(3,4.3)(0,0)
\put(.15,.2){\epsfxsize=3in\epsfbox{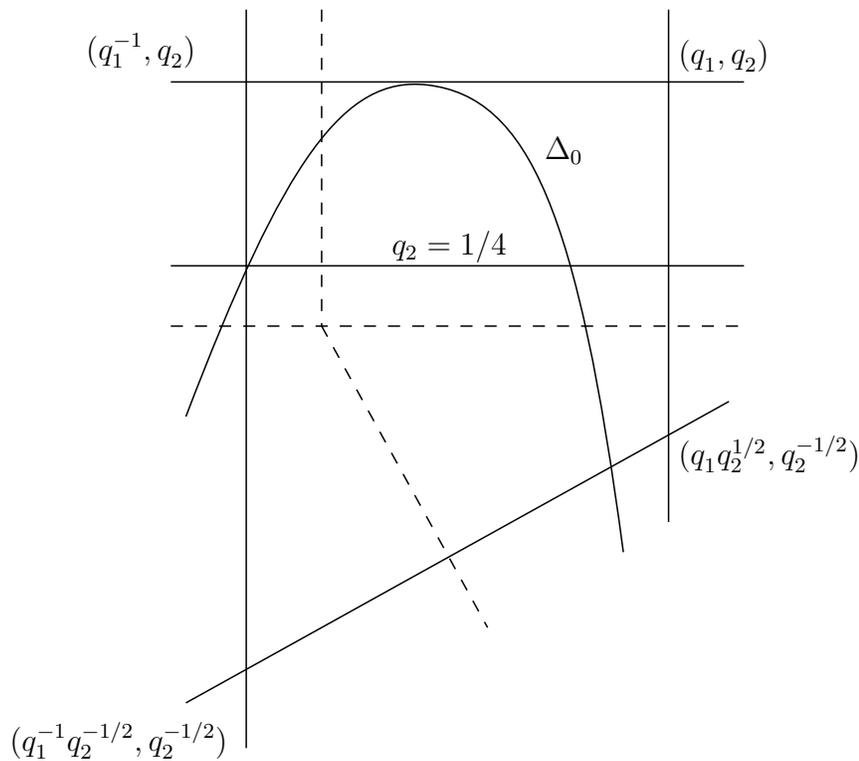}}
\put(2.8,3.8){$(q_1^{\vphantom1},q_2^{\vphantom1})$}
\put(2.8,1.7){$(q_1^{\vphantom1}q_2^{1\smash/2},q_2^{-1\smash/2})$}
\put(-.7,.2){$(q_1^{-1}q_2^{-1\smash/2},q_2^{-1\smash/2})$}
\put(-.3,3.82){$(q_1^{-1},q_2^{\vphantom1})$}
\put(1.3,2.8){$q_2=1/4$}
\put(2.1,3.3){$\Delta_0$}
\end{picture}
\end{center}
\caption{The moduli space for Calabi--Yau hypersurfaces in
$\wp/G$.}\label{fig:1}
\end{figure}
\fi

The moduli space is illustrated in figure \ref{fig:1}.
The figure displays four coordinate charts, which cover the entire moduli
space: the coordinates in each chart are indicated near the point which is
the center of the coordinate chart.  (The dotted lines indicate the
approximate division into ``phase regions,'' described as cones in the
variables ${\frac1{2\pi i}}\log(q_j)$: cf.~\cite{MP}.  Note that our figure
is rotated by $90^\circ$ with respect to figure 5 of \cite{cdfkm}.)  The
chart in
the upper right corner
with coordinates $(q_1,q_2)$ is centered at the so-called ``large complex
structure limit'' point.

The ``discriminant
locus'' where these hypersurfaces become singular has five components: the
principal component, labeled $\Delta_0$ in the figure, is the curve defined by
\begin{equation}
\Delta_0=\left\{ q_2={\frac14}\left(1-{\frac1{256}q_1}\right)^2 \right\};
\end{equation}
the other components are described by $\{q_1=0\}$, $\{q_2=0\}$,
$\{q_1^{-1}=0\}$, and $\{q_2=1/4\}$.\footnote{Note that in \cite{cdfkm},
the principal component was called $C_{\text{con}}$,
and the other components were called $D_{(1,0)}$, $C_\infty$, $C_0$, and $C_1$,
respectively.}  (There is also some monodromy around the ``orbifold locus''
$\{q_2^{-1\smash/2}=0\}$.)

Along the locus where $q_2=1/4$, we found in
\cite{cdfkm} an extremal transition, as follows.
(Notice that this is not a special case of example \ref{subsec:32}, since
we are not simply setting coefficients to zero.)
When $q_2=1/4$,
choose square roots $\sqrt{c_1}$ and $\sqrt{c_2}$ which are related by
requiring that $2\sqrt{c_1}\sqrt{c_2}=c_6$.  Let $\Gamma$ be the image in
$\Aut(\p^5)$ of
\begin{equation*}
\widetilde{\Gamma}:=\{
\vec{\mu}
\in({\Bbb C}^*)^6\ |\
\mu_0^2=\mu_1^4=\mu_2^4=\mu_3^2=\mu_4^2=\mu_5^2
=\mu_1\mu_2=\mu_0\mu_3\mu_4\mu_5=1\}.
\end{equation*}
Then we can
define a rational map $\wp/G\to \p^5/\Gamma$ by
\begin{align*}
y_0&=\sqrt{c_1}\,x_1^4+\sqrt{c_2}\,x_2^4 &y_3&=x_3^2\\
y_1&=x_1\sqrt{x_3x_4x_5} &y_4&=x_4^2\\
y_2&=x_2\sqrt{x_3x_4x_5} &y_5&=x_5^2
\end{align*}
(using the same value of $\sqrt{x_3x_4x_5}$ in both $y_1$ and $y_2$).
The image of $\wp/G$ satisfies the equation
\begin{equation}
\sqrt{c_1}\,y_1^4+\sqrt{c_2}\,y_2^4=y_0y_3y_4y_5,
\end{equation}
while the image of the Calabi--Yau hypersurface additionally
satisfies the equation
\begin{equation}
y_0^2+c_0y_1y_2+c_3y_3^2+c_4y_4^2+c_5y_5^2.
\end{equation}
It can be easily checked that the hypersurface in $\wp/G$ is mapped
{\it birationally}\/ to the complete intersection in $\p^5/\Gamma$, which
is the candidate mirror partner for the complete intersection in $\p^5$ of
bidegree $(2,4)$!

\section{The Location of the Extremal Transition}\label{sec:4}

In order to make the principle formulated in section \ref{sec:2} more precise,
we need to recall the conjectural correspondence between boundary points of
complex structure moduli spaces, and possible mirror partners of a given
Calabi--Yau manifold $X$ \cite{AGM, compact}.
If we compactify the complex structure moduli space
in such a way that the boundary is a divisor with normal crossings, then
candidates for ``large complex structure limit points'' can be identified
by the properties of the monodromy transformations around boundary
divisors.  If the moduli space has dimension $r$, then any such candidate
point $P$ should lie at the intersection of $r$ boundary divisors $D_i$ whose
monodromy transformations $T_1$, \dots, $T_r$ define a ``monodromy weight
filtration'' which is opposite to the Hodge filtration
\cite{deligne, compact}.  The conjecture is that any such point will have an
associated mirror partner $\widehat{Y}$, together with a
choice\footnote{This choice may look a bit unnatural, but it is needed to
get reasonable coordinates; the independence from choices would follow from
the ``cone conjecture.''  See \cite{compact} for a discussion of this issue.}
of a simplicial
rational polyhedral cone $\Pi$ contained in the closure of the
K\"ahler cone of $\widehat{Y}$, in such a way that under the mirror map
$\mu$, the
complement of the $D_i$'s in a neighborhood of $P$ is mapped to an open subset
in the closure $\overline{\cal D}_\Pi$ of the space
\begin{equation}
{\cal D}_\Pi:=
\{\beta\in H^2(\widehat{Y},{\Bbb C}/{\Bbb Z})\ |\ \Im(\beta)\in\Pi\}.
\end{equation}
In the strongest form of the conjecture, one asserts that $\Pi$ can be
chosen so that it is
generated by a basis $e^1$, \dots, $e^r$ of $H^2(\widehat{Y},{\Bbb
Z})/(\text{torsion})$.
In this case, if we write a general element of $H^2(\widehat{Y},{\Bbb C})$
in the form $\sum
t_je^j$ and let $w_j:=\exp(2\pi i\,t_j)$,  we can describe ${\cal D}_\Pi$ in
coordinates as
\begin{equation}
{\cal D}_\Pi=\{(w_1,\dots,w_r)\ |\ 0<|w_j|<1\ \forall j\}.
\end{equation}
This space has a partial compactification
\begin{equation}
{\cal D}_\Pi^-:=\{(w_1,\dots,w_r)\ |\ 0\le|w_j|<1\ \forall j\},
\end{equation}
and $P$ maps to the {\it distinguished limit point}\/
$\mu(P)=(0,\dots,0)\in{\cal D}_\Pi^-$.

If we choose $\Pi$ so that it shares with the K\"ahler cone a face
associated to an extremal transition (say the face spanned by the first
$k$ vectors), then the extremal transition na\"{\i}vely would be expected
at $t_{k+1}=\dots=t_r=0$,
i.e., at $w_{k+1}=\dots=w_r=1$.  The location of the extremal transition
may, however,
be modified by quantum effects in the complexified K\"ahler moduli space.
We expect it to occur at a place where the conformal field theory has a
singularity, and this is measured by poles in the correlation functions of
the quantum theory.

The locus $w_{k+1}=\dots=w_r=1$ will meet the boundary stratum in ${\cal
D}_\Pi^-$ defined by
$w_1=w_2=\dots=w_{k}=0$, and the ``large radius'' approximation to the
quantum field theory should be
good near that boundary stratum.  In fact, we should be able use the
behavior of correlation functions along that boundary stratum
to predict the location of the extremal transition in the complex structure
moduli space.
To do so, we will need to use the ``flat coordinates'' $z_1$, \dots, $z_r$
which are
intrinsically associated to the large complex structure limit point, since
the mirror map $\mu$ has the property that $\mu^*(w_j)=z_j$.

Recall how the flat coordinates are defined: the monodromy properties of
the periods near the large complex structure limit point $P$ guarantee that
if $q_1$, \dots, $q_r$ are local coordinates such that the boundary
divisors intersecting at $P$ are given by $D_j=\{q_j=0\}$, then
there are periods integrals $\varpi_j=\int_{\gamma_j}\Omega$ of the
holomorphic $n$-form $\Omega$ with the property that $\varpi_0$ is
single-valued near $P$, while
\begin{equation}
\varpi_j=\frac{\varpi_0}{2\pi i}\log(q_j)+ \text{single-valued function}.
\end{equation}
The flat coordinates are then given by $z_j=\exp(2\pi
i{\varpi_j}/{\varpi_0})$; they are uniquely determined up to
multiplication by constants.\footnote{For a discussion of how those
constants should be fixed, see \cite{predictions}.}

To learn what we should expect concerning the location of the extremal
transition, let us again consider the example from
section \ref{subsec:33}.  In that example, we should study the boundary
curve $B$ (in the complex structure moduli space of $X$) defined by
$z_1=0$, or equivalently by $q_1=0$.  The period integrals of $X$ satisfy
certain Gelfand-Kapranov-Zelevinsky hypergeometric
differential equations \cite{batyrev:vmhs}; when restricted to the locus
$q_1=0$, there is a single such equation, which can be read off of the
formula $q_2=c_1c_2/c_6^2$ as being
\begin{equation}
\left((q_2\frac{d}{dq_2})(q_2\frac{d}{dq_2}) -
q_2(-2q_2\frac{d}{dq_2})(-2q_2\frac{d}{dq_2}-1)\right)\varpi(0,q_2)=0.
\end{equation}
This has a general solution near $q_2=0$ given by \cite{Small}
\begin{equation}
\varpi(0,q_2)
=C_1+C_2\log\left(\frac{2q_2}{1-2q_2+\sqrt{1-4q_2}}\right),
\end{equation}
choosing the branch of the square root which is near $1$ when $q_2$ is near
$0$. It follow that the flat coordinate along $B$ is given by
\begin{equation}
z_2|_B=\frac{2q_2}{1-2q_2+\sqrt{1-4q_2}}.
\end{equation}
Note that (as can be seen in figure \ref{fig:1}), $B$ meets two other
components
of the discriminant locus, at $q_2=1/4$ and at $q_2\to\infty$ (the latter
being the intersection with the ``orbifold locus'').  When
$q_2=1/4$, we have $z_2=1$ whereas when $q_2\to\infty$ we have $z_2=-1$.
Thus, there will be poles in correlation functions precisely at $z_2=\pm1$.

More generally, we should expect that poles in correlation functions could
occur at several distinct values of $|z_r|$.  (This is known to happen in
other examples \cite{Small, gmv}.)  The large
radius approximation can only be trusted for values of $|w_r|$ less than
the minimum value at which a pole occurs, so we shall expect that any
extremal transition whose occurrence is predicted by mirror symmetry will
occur at a pole where the value of $|z_r|=|\mu^*(w_r)|$ is minimal.  As the
present example shows, such a pole need not be unique.

Returning to the general case, we formulate the following conjecture
concerning the location of the extremal transition, which we hope is not
too far off the mark.  We assume that a mirror pair $(X,\widehat{Y})$ is
somehow known, corresponding to the large complex structure limit point $P$
(for $X$) and the subcone $\Pi$ of the K\"ahler cone of $\widehat{Y}$.  Let
${\cal F}^\perp=\Pi\cap\Span\{e^1,\dots,e^k\}$ be a face of $\Pi$.

\begin{conjecture}
There exist a compactification $\overline{\cal M}$ of the complex structure
moduli space of $X$ containing $P=D_1\cap\dots\cap D_r$, together with
components $\Delta_{k+1}$, \dots, $\Delta_r$ of the boundary of
$\overline{\cal M}$ along which $X$ acquires canonical singularities, and a
stratification
\begin{equation}
\Delta_{k+1}\cup\dots\cup\Delta_r=\coprod_{\sigma<{\cal F}} \Delta_\sigma
\end{equation}
of the union of those components, indexed by subcones of the dual face
${\cal F}$ of the Mori
cone, such that the intersection of the stratum $\Delta_\sigma$ with
\begin{equation}
B_\sigma:=D_{\sigma(1)}\cap\dots\cap
D_{\sigma(s)}=\{z_{\sigma(1)}=\dots=z_{\sigma(s)}\}
\end{equation}
lies at the minimum possible distance from the origin, among locations of
poles of correlation functions on $B_\sigma$.

Furthermore, the singular Calabi--Yau space $\overline{X}_{\Delta_{\cal
F}}$ has a Calabi--Yau desingularization
$\widehat{X}\to\overline{X}_{\Delta_{\cal F}}$ if and only if the
contracted space $\overline{Y}_{\cal F}$ has a Calabi--Yau smoothing $Y$.
In this case, the Calabi--Yau manifolds $\widehat{X}$ and $Y$ should be
mirror partners.
\end{conjecture}

We have limited our discussion to neighborhoods of large complex structure
limit points, which are the mirrors of K\"ahler cones (of various
birational models of $\widehat{Y}$).  It is frequently
possible to analytically continue the complexified K\"ahler moduli space
beyond these K\"ahler cones \cite{beyond}, but we don't
know good criteria for deciding about the existence of extremal transitions
(on the ``K\"ahler moduli'' side) in such regions of the moduli space.
Extremal transitions in such regions, if they exist, would evade detection
in the sort of analysis given here.

\section{Recent developments}\label{subsec:recent}

In the physics literature, conifold transitions were first observed in a
process known as ``splitting'' which related various families of complete
intersection Calabi--Yau threefolds \cite{cdls}.  Considerable effort was
expended in showing that all then-known examples of Calabi--Yau threefolds
could be connected into a single web
\cite{ghone, ghtwo, cgh}, and it was also observed that these connections
occurred at finite distance in the moduli space (with respect to the
natural ``Zamolodchikov'' metric on that space) \cite{cghfinite}.  However, a
physical mechanism implementing these transitions was unknown.

Some months after this lecture was delivered, a new mechanism was
proposed in the physics literature for realizing an extremal transition
as a physical process in type II string theory \cite{strom, gms}.
The mechanism in its original form only applies to conifold transitions,
but there are now indications \cite{KMP, MVII, Wittennewest, morsei, GMS,
IMS} that similar mechanisms will enable all extremal transitions to be
realized in physics.  Motivated by this, the subpolyhedron construction
described in section \ref{subsec:32} was subsequently used \cite{CGGK,
ACJM} to show that all Calabi--Yau hypersurfaces in weighted projective
spaces can be linked by extremal transitions.

In another direction---perhaps closer in spirit to the original approach of
Clemens and Reid---Kontsevich \cite{kontsevich} has made a fascinating
construction involving Lagrangian analytic cones in an infinite-dimensional
space, and has conjectured a mirror symmetry relationship in terms of these
cones which would involve {\it all}\/
symplectic complex threefolds with trivial canonical bundle (even those for
which the symplectic structure does not arise from a K\"ahler structure).

Finally, there was been a recent geometric reformulation of the basic
mirror symmetry property in physics \cite{SYZ}
(see also \cite{underlying, GW}), in
terms of fibrations of a Calabi--Yau manifold by special Lagrangian tori.
It is an important and challenging problem to understand how such
fibrations behave under an extremal transition.  Such an understanding
could ultimately lead to a proof of the conjecture in section \ref{sec:4}
using the new ``geometric'' definition of mirror symmetry.

\subsection*{ Acknowledgments}\label{subsec:ack}

I am grateful to Mark Gross for insisting that the analysis given
here should not be limited to the case of small resolutions, and
for discussions on other points.
I would also like to thank Paul Aspinwall, Bob Friedman,
Sheldon Katz, Ronen Plesser, Miles Reid, Pelham Wilson,
Edward Witten
and especially Brian Greene
for useful discussions.

\def\bysame{\leavevmode\hbox to3em{\hrulefill}\thinspace}

\affil{Department of Mathematics, Duke University, Durham, NC 27708-0320 USA}

\end{document}